\input amstex
\documentstyle{amsppt}
\pagewidth{6.0in}\vsize8.5in\parindent=6mm
\parskip=3pt\baselineskip=14pt\tolerance=10000\hbadness=500
\document
\topmatter
\title Note on the Dirichlet Approximation Theorem
\endtitle
\author Yong-Cheol Kim \endauthor
\abstract We survey the classical results on the Dirichlet
Approximation Theorem
\endabstract
\endtopmatter

\define\inn#1#2{\langle#1,#2\rangle}

\define\lcontr{\rfloor}
\define\lco#1#2{{#1}\lcontr{#2}}
\define\lcoi#1#2{\imath({#1}){#2}}
\define\rco#1#2{{#1}\rcontr{#2}}
 \redefine\mat{{\text{\rm Mat}}}
 \redefine\mod{{\text{\rm mod}}}

\define\ap{\alpha}             
\define\bt{\beta}
             
\define\dt{\delta}             
\define\vep{\varepsilon}
\define\zt{\zeta}

\define\kp{\kappa}

\define\sm{\sigma}

            \define\iy{\infty}
\define\lt{\left}            \define\rt{\right}
\define\f{\frac}             \define\el{\ell}


\define\BN{{\Bbb N}}

\define\BQ{{\Bbb Q}}
\define\BR{{\Bbb R}}

\define\BZ{{\Bbb Z}}

\define\cF{{\Cal F}}




\define\s{\setminus}         
            \define\e{\eta}
        \define\fd{\fallingdotseq}
       
     \define\ds{\dsize} 

In this chapter, we show up the theory of Diophantine
approximation leads to understand how well real numbers can be
captured by relations with the integers. In what follows, we
denote by $\BN_0\fd\BN\cup\{0\}$.

\proclaim{Theorem 1.1 [Dirichlet's Approximation Theorem]}

\noindent For each $\ap\in\BR$ and $N\in\BN$, there are
$n\in\BN\,\,( n\le N )$ and $p\in\BZ$ such that
$$\lt|\ap-\f{p}{n}\rt|<\f{1}{N n},\,\,\,\,i.e.\,\,\, |n\ap-p|<\f{1}{N}.$$
\endproclaim

\noindent{\it Proof.} For $n\in\BN_0$, we set $x_n=n\ap-[n\ap]$.
We also divide $[0,1)$ into $N$ half open subintervals of equal
length $$I_m=\lt[\f{m-1}{N},\f{m}{N}\rt),\,\,m=1,2,\cdots,N.$$
Since $\{x_0,x_1,\cdots,x_N\}\subset\bigcup_{m=1}^N I_m$, at least
two of $x_n$'s belongs to a common interval $I_{m_0}$. Let $x_k,
x_{\el}$ be such two numbers inside $I_{m_0}$ and let $k<\el$.
Then we have that $$\f{m_0-1}{N}\le
k\ap-[k\ap]<\f{m_0}{N}\,\,\text{ and
}\,\,\f{m_0-1}{N}\le\el\ap-[\el\ap]<\f{m_0}{N}.$$ By substracting
those two inequalties, we obtain that
$$-\f{1}{N}<\el\ap-k\ap-([\el\ap]-[k\ap])<\f{1}{N}.$$ Then set
$n=\el-k\le N$ and $p=[\el\ap]-[k\ap]$. \qed

\proclaim{Corollary 1.2 [Main theorem on the linear Diophantine
equation]}

\noindent For each $\ds\f{a}{b}\in\BQ$ in lowest terms, there are
$x,y\in\BZ$ such that $a x-b y=1$.
\endproclaim

\noindent{\it Proof.} If $b=1$, then $ax-y=1$ is solved by setting
$x=0$ and $y=-1$. Without loss of generality, we may assume that
$b\ge 2$. Applying Theorem 1.1 with $N=b-1$, there are
$n\in\BN\,\,(n\le N)$ and $p\in\BN$ such that $$|\ap
n-p|=\lt|\f{a}{b} n-p\rt|<\f{1}{N}=\f{1}{b-1}.$$ Multiplying it by
$b$ leads to $\ds |a n-b p|<\f{b}{b-1}=1+\f{1}{b-1}\le 2$. Since
$an-bp\in\BZ$, this implies that $|an-bp|\le 1$.

The case that $an-bp=0$ is excluded, since it implies that
$\ap=\ds\f{a}{b}=\f{p}{n}$ and $n\le N=b-1<b$, which contradicts
to the choice of $\ds\f{a}{b}$. Thus the only possibility is
$an-bp=\pm 1$.

In case that $an-bp=1$, we set $x=n$ and $y=p$; in case that
$an-bp=-1$, we set $x=-n$ and $y=-p$. \qed

\proclaim{Corollary 1.3} The special Pell equation $x^2 -c y^2 =1$
can be solved in the cases $c<0$ and $c=a^2$ for $a\in\BZ$.
\endproclaim

\noindent{\it Proof.} In case that $c=-1$ and $x^2 +y^2 =1$, it
can be solved by $x=\pm 1,\, y=0$ or $x=0,\, y=\pm 1$. In case
that $c<-1$ and $x^2 + |c| y^2=1$, it can be solved by $x=\pm 1,\,
y=0$. Finally, in case that $c=a^2$ and $x^2 -c y^2
=(x-ay)(x+ay)=1$, it can be solved by $x=\pm 1,\, y=0$ for $a\neq
0$, or $x=\pm 1,\, y=t$ with $t\in\BZ$ for $a=0$. \qed

\proclaim{Lemma 1.4} For each $\ap\in\BQ^c$, there are infinitely
many rationals $\ds\f{p}{n}$ in lowest terms such that
$$\lt|\ap-\f{p}{n}\rt|<\f{1}{n^2}.$$
\endproclaim

\noindent{\it Proof.} Assume that there is a $\ap_0\in\BQ^c$ such
that $$0<\lt|\ap_0 -\f{p_k}{n_k}\rt|<\f{1}{n_k^2}\,\,\text{ for
all }k=1,2,\cdots,K.$$ Then we set $\ds\dt=\min_{1\le k\le
K}\lt[\,n_k\lt|\ap_0 -\f{p_k}{n_k}\rt|\,\,\rt]$. Applying Theorem
1.1 with $N=\ds\lt[\f{1}{\dt}\rt]+1>\f{1}{\dt}$, there are
$n\in\BN\,\,(n\le N)$ and $p\in\BZ$ such that $$\lt|\ap_0
-\f{p}{n}\rt|<\f{1}{nN}.$$ Without loss of generality,
$\ds\f{p}{n}$ may here assumed to be in lowest terms. Then we have
that $$\lt|\ap_0 -\f{p}{n}\rt|\le\f{1}{n^2},$$ and so
$\ds\f{p}{n}$ is one of $\ds\f{p_k}{n_k}$'s. Thus, $p=p_k$ and
$n=n_k$ for some $k, 1\le k\le N$, and thus we obtain that
$$n_k\lt|\ap_0 -\f{p_k}{n_k}\rt|<\f{1}{N}<\dt;$$ which contradicts
to the definition of $\dt$. \qed

\proclaim{Corollary 1.5 [Fermat's theorem on the Pell equation]}

\noindent If $c\in\BN$ and $c\neq a^2$ for any $a\in\BZ$, then the
special Pell equation $x^2 -c y^2=1$ has infinitely many integral
solutions.
\endproclaim

\noindent{\it Proof.} First of all, we show that there is a single
solution $x=\xi,\, y=\e$ with $\e\neq 0$. By Lemma 1.4, there are
infinitely many rationals $\ds\f{p}{n}$ in lowest terms such that
$$0<\lt|\sqrt c -\f{p}{n}\rt|<\f{1}{n^2}.$$ If we set
$\ap=p+n\sqrt c$ and $\overline \ap=p-n\sqrt c$, then we have $\ds
|\overline\ap|<\f{1}{n}$, and so $$|\ap|=|p-n\sqrt c+2n\sqrt
c\,|\le |\overline\ap|+2n\sqrt c\le\f{1}{n}+2n\sqrt c.$$ Thus we
have that $\ds |p^2 -c
n^2|=|\ap\overline\ap|\le\lt(\f{1}{n}+2n\sqrt
c\rt)\f{1}{n}=\f{1}{n^2}+2\sqrt c\le 2\sqrt c +1.$ Since the
infinitely many integers $p^2 -c n^2$ lie in $[-2\sqrt c -1,
2\sqrt c +1]$, there is $a\in\BZ\s\{0\}$ which coincides with
infinitely many of the $p^2 -c n^2$. Set $\cF=\{p+n\sqrt c : a=p^2
-c n^2\}$. Then we define the equivalence relation on $\cF$ by
$$\ap=p+n\sqrt c\,\sim\,\bt=q+m\sqrt c \,\,\Leftrightarrow\,\, p\equiv
q,\,n\equiv m \,\,(\,\mod \,|a|\,).$$ There are at most $a^2$
equivalence classes in the quotient space $\cF/\sim$, and so in at
least one of these finitely many equivalence classes there must be
infinitely many of the numbers $\ap,\bt,\cdots$ as above. Take
such an equivalent pair $\ap=p+n\sqrt c,\,\bt=q+m\sqrt c$ with
$\ap\neq\bt$; say, $|\ap|>|\bt|$. So the number
$$\f{\ap-\bt}{a}=\f{p-q}{a}+\f{n-m}{a}\sqrt c$$ is of the form
$x+y\sqrt c,\,x,y\in\BZ$. Thus the number
$$\vep\fd\f{\ap}{\bt}=1+\f{a}{b}\cdot\f{\ap-\bt}{a}=1+\f{\bt\overline\bt}{\bt}\cdot\f{\ap-\bt}{a}=
1+(q-m\sqrt c)\lt(\f{p-q}{a}+\f{n-m}{a}\sqrt c\rt)$$ is of the
form $\vep=\xi+\e\sqrt c$ with $\xi,\e\in\BZ$, for which
$$\xi^2-c\e^2=\vep\overline\vep=\f{\ap\overline\ap}{\bt\overline\bt}=\f{a}{a}=1.$$
Since $|\vep|>1$, we see that $\e\neq 0$, and thus $\xi,\,\e$ is
the required solution.

For the remaining part, we write $\vep=\xi+\e\sqrt c$ and
$\overline\vep=\xi-\e\sqrt c$. Then we have that $$1=\xi^2-c
\e^2=(\xi+\e\sqrt c)(\xi-\e\sqrt c)=\vep\cdot\overline\vep,$$ and
so $\vep^n\cdot{\overline\vep}^n=1$ for all $n\in\BN$. Moreover,
the number $\vep^n=(\xi+\e\sqrt c)^n$ is of the form $x+y\sqrt c$
with $x,y\in\BZ$, which is represented uniquely by the
irrationality of $\sqrt c$, and also $\vep^n=x+y\sqrt c$ and
${\overline\vep}^n=x-y\sqrt c$. Hence the infinitely many powers
$\vep^n=x+y\sqrt c$ lead to infinitely many solutions $x,y$ of the
special Pell equation, which are of course all distinct from each
other since $|\vep|\neq 1$. \qed

\proclaim{Proposition 1.6} $\ap\in\BQ^c$ if and only if for each
$\vep>0$, there are $x,y\in\BZ$ such that $0<|\ap x-y|<\vep$.
\endproclaim

\noindent{\it Proof.} $(\Rightarrow)$ By Lemma 1.4, there are
infinitely many rationals $\ds\f{p_n}{q_n}$ with strictly
increasing denominators satisfying
$\ds\lt|\ap-\f{p_n}{q_n}\rt|<\f{1}{q_n^2}$. Given $\vep>0$, we
choose $n\in\BN$ so large that $\ds q_n>\f{1}{\vep}$. Set $x=q_n$
and $y=p_n$. Then we have that $0<|\ap x-y|<\vep$.

$(\Leftarrow)$ Assume that the righthand side holds. If
$\ap\in\BQ$, i.e. $\ap=\ds\f{a}{b}$ with $a\in\BZ$ and $b\in\BN$,
then we have that $$0<\lt|\f{a}{b} x-y\rt|<\f{1}{b},$$ i.e.
$0<|ax-by|<1$, which is impossible. \qed

\proclaim{Example} Set $\ap$ equal to the Cantor series
$$\ap=\sum_{n=1}^{\iy}\f{z_n}{g_1 g_2\cdots
g_n}=\f{z_1}{g_1}+\f{z_2}{g_1 g_2}+\f{z_3}{g_1 g_2 g_3}+\cdots$$
with an increasing sequence of natural numbers $2\le g_1\le
g_2\le\cdots\le g_n\le\cdots$, $z_n\in\{0,1\}$, and $z_n=1$ for
infinitely many $n$. Then $\ap\in\BQ^c$ and the cardinal number of
such $\ap$'s is uncountable.
\endproclaim

\noindent{\it Proof.} For $N\in\BN$, we set $G_N=g_1 g_2\cdots
g_N$ and $\ds\sum_{n=1}^N\f{z_n}{g_1 g_2\cdots g_n}=\f{P_N}{G_N}$.
Then we have that $$ \split
0<\lt|\ap-\f{P_N}{Q_N}\rt|&=\lt|\f{z_{N+1}}{g_1\cdots g_N\,
g_{N+1}}+\f{z_{N+2}}{g_1\cdots g_N\, g_{N+1}\, g_{N+2}}+\cdots \rt| \\
&\le\f{1}{g_1\cdots g_N\,
g_{N+1}}\lt(1+\f{1}{g_{N+2}}+\f{1}{g_{N+2}\,
g_{N+3}}+\cdots \rt) \\
&\le\f{1}{G_N\, g_{N+1}}\lt(
1+\f{1}{g_{N+1}}+\f{1}{g^2_{N+1}}+\cdots \rt) \\
&=\f{1}{G_N
g_{N+1}}\cdot\f{1}{1-\ds\f{1}{g_{N+1}}}=\f{1}{G_N(g_{N+1}-1)}.
\endsplit $$ If the sequence $\{g_n\}$ is unbounded, then $$0<|\ap
G_N-P_N|\le\f{1}{g_{N+1}-1}$$ can be made to be arbitrarily small
for all large enough $N$. So it follows from Proposition 1.6. \qed

\proclaim{Corollary 1.7} $e\in\BQ^c$. \endproclaim

\noindent{\it Proof.} It easily follows from taking $g_n=n$ and
$z_n=1$ for all $n\in\BN$ in the above example. \qed

Next we show up the irrationality of the second fundamental
constant associated with the circle, i.e. $\pi$,
$\zt(2)=\ds\sum_{n=1}^{\iy}\f{1}{n^2}=\f{1}{6}\pi^2$, and
$\zt(3)=\ds\sum_{n=1}^{\iy}\f{1}{n^3}$. The essential core of the
proof is derived from the following formulae $$ \split
\zt(2)&=1+\f{1}{2^2}+\cdots
+\f{1}{r^2}+\int_0^1\int_0^1\f{(xy)^r}{1-xy}\,dxdy, \\
\zt(3)&=1+\f{1}{2^3}+\cdots
+\f{1}{r^3}-\f{1}{2}\int_0^1\int_0^1\f{\ln(xy)}{1-xy}(xy)^r\,dxdy.
\endsplit $$
We obtain that for $r,s,\sm\in\BN_0$,
$$\int_0^1\int_0^1\f{x^{r+\sm}y^{s+\sm}}{1-xy}\,dx\,dy=\sum_{k=0}^{\iy}
\f{1}{(k+r+\sm+1)(k+s+\sm+1)}\tag{1.1} $$ by expanding
$(1-xy)^{-1}$ as a geometric series. If $r=s$ and $\sm=0$ in
(1.1), then we have that
$$\int_0^1\int_0^1\f{(xy)^r}{1-xy}\,dx\,dy=\zt(2)-\lt(1+\f{1}{2^2}+\cdots
+\f{1}{n^2}\rt).\tag{1.2}$$ If $r>s$ in (1.1), then we have that
$$ \split \int_0^1\int_0^1\f{x^{r+\sm}
y^{s+\sm}}{1-xy}\,dx\,dy&=\sum_{k=0}^{\iy}\lt(\f{1}{k+s+\sm+1}-\f{1}{k+r+\sm+1}\rt)
\\ &=\f{1}{r-s}\lt(\f{1}{s+1+\sm}+\cdots
+\f{1}{r+\sm}\rt), \endsplit \tag{1.3} $$ which represents a
rational number whose denominator is certainly given by
${V(r)}^2$, the square of the least common multiple $V(r)$ of
$1,2,\cdots, r$. Now we differentiate the both sides of (1.3) with
respect to $\sm$, and then set $\sm=0$. Then we obtain that
$$\int_0^1\int_0^1\f{\ln(xy)}{1-xy}\, x^r
y^s\,dx\,dy=\f{-1}{r-s}\lt(\f{1}{(s+1)^2}+\cdots +\f{1}{r^2}\rt),
\,r>s,$$ which describes a rational number whose denominator is
certainly given by ${V(r)}^3$. We differentiate the both sides of
(1.1) on the case $r=s$ with respect to $\sm$, and then set
$\sm=0$. Then we have that
$$ \split
\int_0^1\int_0^1\f{\ln(xy)}{1-xy}\,(xy)^r\,dx\,dy&=\sum_{k=0}^{\iy}\f{-2}{(k+r+1)^3}
\\ &=-2\lt[\zt(3)-\lt(1+\f{1}{2^3}+\cdots +\f{1}{r^3}\rt)\rt]
\endsplit \tag{1.4} $$ So it easily follows from (1.2) and (1.3)
that $$\int_0^1\int_0^1\f{(1-y)^n P_n(x)}{1-xy}\,dx\,dy=\f{a_n
\zt(2)+b_n}{{V(n)}^2}$$ where
$P_n(x)=\ds\f{1}{n!}\lt(\f{d}{dx}\rt)^n (x^n (1-x)^n )$ and $a_n,
b_n$ are some integers. Here we observe that $P_n(x)$ has only
integral coefficients. Moreover, the integration by parts
$n$-times with respect to $x$ leads us to get
$$\f{a_n\zt(2)+b_n}{{V(n)}^2}=(-1)^n\int_0^1\int_0^1\f{y^n (1-y)^n
x^n (1-x)^n}{(1-xy)^{n+1}}\,dx\,dy.$$ Calculating the maximum
value of the function $\ds\f{y(1-y)x(1-x)}{1-xy}$ on $[0,1]\times
[0,1]$, we obtain that $$0\le\f{y(1-y)x(1-x)}{1-xy}\le\lt(\f{\sqrt
5 -1}{2}\rt)^5,$$ which implies that $$ \split
\f{|a_n\zt(2)+b_n|}{{V(n)}^2}&=\lt|\int_0^1\int_0^1\f{y^n (1-y)^n
x^n (1-x)^n}{(1-xy)^{n+1}}\,dx\,dy\rt| \\
&\le\lt(\f{\sqrt 5
-1}{2}\rt)^{5n}\int_0^1\int_0^1\f{1}{1-xy}\,dx\,dy\le\lt(\f{\sqrt
5 -1}{2}\rt)^{5n}\zt(2). \endsplit $$ Thus we obtain that
$$0<|a_n\zt(2)+b_n|\le {V(n)}^2\lt(\f{\sqrt 5
-1}{2}\rt)^{5n}\zt(2).\tag{1.5}$$ If we can show that
$\ds\lim_{n\to\iy}{V(n)}^2\lt(\f{\sqrt 5 -1}{2}\rt)^{5n}=0$, we
can apply Proposition 1.6 to conclude it because the righthand
side of (1.5) could be chosen to be arbitrarily small. For this,
we observe that $\ds V(n)\le n^{\pi(n)}$ where $\pi(n)$ denotes
the number of primes less than or equal to $n$. By the prime
number theorem to be shown in Chapter 5, we see that
$$\lim_{n\to\iy}\f{\pi(n)}{n/\ln n}=1.$$ Thus we have that
$\ds\pi(n)\le\ln 3\cdot\f{n}{\ln n}$, and so $\ds V(n)\le
n^{n^{\ln 3/\ln n}}=3^n$ for all sufficiently large $n$. Hence we
conclude that $$0\le\lim_{n\to\iy}{V(n)}^2\lt(\f{\sqrt 5
-1}{2}\rt)^{5n}\le\lim_{n\to\iy} 9^n\lt(\f{\sqrt 5
-1}{2}\rt)^{5n}\le\lim_{n\to\iy}\lt(\f{5}{6}\rt)^n=0.$$ This
proves Corollary 1.8 and Corollary 1.9.

\proclaim{Corollary 1.8} $\zt(2)=\ds\f{1}{6}\pi^2\in\BQ^c$.
\endproclaim

\proclaim{Corollary 1.9} $\pi\in\BQ^c.$ \endproclaim

Let us check up the irrationality of $\zt(3)$. As in the above
process, we obtain that
$$\int_0^1\int_0^1\f{-\ln(xy)}{1-xy}\,P_n(x)\,P_n(y)\,dx\,dy=\f{a_n\zt(3)+b_n}{{V(n)}^3}$$
where $a_n$ and $b_n$ are some integers. The identity
$\ds\f{-\ln(xy)}{1-xy}=\int_0^1\f{1}{1-(1-xy)z}\,dz$ leads to get
that $$ \split &\f{a_n\zt(3)+b_n}{{V(n)}^3} \\
&=\int_0^1\int_0^1\int_0^1\f{P_n(x)\,P_n(y)}{1-(1-xy)z}\,dx\,dy\,dz
\\ &=\int_0^1\int_0^1\int_0^1\f{(xyz)^n (1-x)^n\,P_n(y)}{(1-(1-xy)z)^{n+1}}\,dx\,dy\,dz
\,\,\text{( integration by parts $n$-times w.r.t. $x$ )}\\
&=\int_0^1\int_0^1\int_0^1 (1-x)^n (1-w)^n
\f{P_n(y)}{1-(1-xy)w}\,dx\,dy\,dw
\,\,\text{( change of variable $w=\f{1-z}{1-(1-xy)z}$ )}\\
&=\int_0^1\int_0^1\int_0^1\f{x^n(1-x)^n y^n(1-y)^n w^n(1-w)^n
}{(1-(1-xy)w)^{n+1}}\,dx\,dy\,dw \,\,\text{( integration by parts
$n$-times w.r.t. $y$ ).} \endsplit $$ Computing the extreme value
of the function $\ds\f{x(1-x) y(1-y) w(1-w)}{1-(1-xy)w}$ on
$[0,1]^3$, we have that $$\f{x(1-x) y(1-y) w(1-w)}{1-(1-xy)w}\le
(\sqrt 2 -1)^4.$$ Thus it follows from this and (1.4) that $$
\split 0\le\f{|a_n\zt(3)+b_n|}{{V(n)}^3}&\le (\sqrt 2 -1)^{4n}
\int_0^1\int_0^1\int_0^1\f{1}{1-(1-xy)w}\,dx\,dy\,dw \\
&=(\sqrt 2 -1)^{4n}\int_0^1\int_0^1\f{-\ln(xy)}{1-xy}\,dx\,dy \\
&=2\,\zt(3)\cdot (\sqrt 2-1)^{4n}. \endsplit $$ Since $V(n)\le
3^n$ for all sufficiently large $n$, we have that
$$0\le |a_n\zt(3)+b_n|\le 2\,\zt(3)\cdot {27}^n\cdot (\sqrt
2-1)^{4n}<2\,\zt(3)\lt(\f{4}{5}\rt)^n\rightarrow 0\,\,\text{ as
$n\to\iy$. }$$ Hence this proves Corollary 1.10.

\proclaim{Corollary 1.10 [Ap\'ery's Theorem]}
$\ds\zt(3)=1+\f{1}{2^3}+\f{1}{3^3}+\f{1}{4^3}+\cdots \in\BQ^c.$
\endproclaim

We next generalize the Dirichlet Approximation Theorem in
Proposition 1.11 and Proposition 1.12. Moreover we prove its
generalized formation containing them in Theorem 1.13 which is due
to Kronecker.

\proclaim{Proposition 1.11} For each
$(\ap_1,\cdots,\ap_L)\in\BR^L$ and $N\in\BN$, there are
$n\in\BN\,\,(n\le N^L)$ and $(p_1,\cdots,p_L)\in\BZ^L$ such that
$$\max_{1\le\el\le L}\lt|\ap_{\el}-\f{p_{\el}}{n}\rt|<\f{1}{Nn}.$$
\endproclaim

\proclaim{Proposition 1.12 [due to Dirichlet]} For each
$(\ap_1,\cdots,\ap_L)\in\BR^L$ and $N\in\BN$, there are $p\in\BN$
and $(n_1,\cdots,n_L)\in\BZ^L\s\{0\}$ with $\ds\max_{1\le\el\le
L}|n_{\el}|\le N^{1/L}$ such that $$\lt|\,\sum_{\el=1}^L
\ap_{\el}n_{\el}-p\,\rt|<\f{1}{N}.$$
\endproclaim

\proclaim{Theorem 1.13 [Multidimensional Dirichlet Approximation
Theoerm]}

\noindent For each
$(\ap_{11},\ap_{12},\cdots,\ap_{ML})\in\BR^{ML}$ and $N\in\BN$,
there are $(p_1,\cdots,p_M)\in\BZ^M$ and
$(n_1,\cdots,n_L)\in\BZ^L\s\{0\}$ with $\ds\max_{1\le\el\le
L}|n_{\el}|\le N^{M/L}$ such that $$\max_{1\le m\le
M}\lt|\,\sum_{\el=1}^L\ap_{m\el} n_{\el}-p_m\,\rt|<\f{1}{N}.$$
\endproclaim

\noindent{\it Proof.} We observe that
$$\ds\lt(\,\,\,\sum_{\el=1}^L\ap_{1\el}\,n_{\el}-\lt[\sum_{\el=1}^L\ap_{1\el}\,n_{\el}\rt],
\cdots,\sum_{\el=1}^L\ap_{M\el}\,n_{\el}-\lt[\sum_{\el=1}^L\ap_{M\el}\,n_{\el}\rt]\,\,\,\rt)\in
[0,1)^M.\tag{1.6}$$ We decompose $[0,1)^M$ into $N^M$ subcubes
$$Q_{k_1,\cdots,k_M}=\lt[\f{k_1
-1}{N},\f{k_1}{N}\rt)\times\cdots\times\lt[\f{k_M
-1}{N},\f{k_M}{N}\rt),\,\,1\le k_m\le N.$$ If we fix $0\le
n_{\el}\le P$, then there are $(P+1)^L$ $L$-tuples
$(n_1,\cdots,n_L)$; moreover, the requirement for applying the
Dirichlet pigeon hole principle is $(P+1)^L >N^M$, i.e. $P+1>
N^{M/L}$.  Set $P=[N^{M/L}]$. By the pigeon hole principle, two of
such $M$-tuples in (1.6) lie in one of the subcubes
$Q_{k_1,\cdots,k_M}$; i.e. for all $m=1,2,\cdots,M$ and for
distinct $n'_{\el}$ and $n''_{\el}$,
$$ \split
\f{k_m-1}{N}&\le\sum_{\el=1}^L\ap_{m\el}\,n'_{\el}-\lt[\sum_{\el=1}^L\ap_{m\el}\,n'_{\el}\rt]
<\f{k_m}{N}, \\
\f{k_m-1}{N}&\le\sum_{\el=1}^L\ap_{m\el}\,n''_{\el}-\lt[\sum_{\el=1}^L\ap_{m\el}\,n''_{\el}\rt]
<\f{k_m}{N}. \endsplit $$ Set $n_{\el}=n''_{\el}-n'_{\el}$ and
$p_m=\ds\lt[\sum_{\el=1}^L\ap_{m\el}\,n''_{\el}\rt]-\lt[\sum_{\el=1}^L\ap_{m\el}\,n'_{\el}\rt].$
Then we obtain that $$\max_{1\le m\le
M}\lt|\,\sum_{\el=1}^L\ap_{m\el} n_{\el}-p_m\,\rt|<\f{1}{N}.$$
Since $n'_{\el}, n''_{\el}$ are distinct and $0\le
n'_{\el},n''_{\el}\le P\le N^{M/L}$, not all the $n_{\el}$ are
equal to zero and $\ds\max_{1\le\el\le L}|n_{\el}|\le N^{M/L}$.
\qed

In the following corollary, one can ask how small linear forms
$\ds\lt(\,\,\sum_{\el=1}^L\ap_{1\el}\,
x_{\el},\cdots,\sum_{\el=1}^L\ap_{M\el}\,x_{\el}\,\,\rt)$ can
occur.

\proclaim{Corollary 1.14} Foe each $N\in\BN$ and each $M$-tuples
$\ds\lt(\,\,\sum_{\el=1}^L\ap_{1\el}\,
x_{\el},\cdots,\sum_{\el=1}^L\ap_{M\el}\,x_{\el}\,\,\rt)$ of
linear forms derived from $(\ap_{m\el})\in\mat_{M\times
L}(\BR)\,(M<L)$, there is $(x_1,\cdots,x_L)\in\BZ^L\s\{0\}$ with
$\ds\max_{1\le\el\le L}|x_{\el}|\le N^{M/L}$ such that
$$\max_{1\le m\le M}\lt|\sum_{\el=1}^L\ap_{m\el}\,x_{\el}\rt|\le
2\,L\cdot A\cdot N^{M/L -1}$$ where $A=\ds\max_{1\le m\le
M}|\ap_{m\el}|$.
\endproclaim

\noindent{\it Proof.} By Theorem 1.13, there are
$(p_1,\cdots,p_M)\in\BZ^M$ and $(x_1,\cdots,x_L)\in\BZ^L\s\{0\}$
with $\ds\max_{1\le\el\le L}|x_{\el}|\le N^{M/L}$ such that
$$\max_{1\le m\le M}\lt|\sum_{\el=1}^L\kp\,\ap_{m\el}\,
x_{\el}-p_m\rt|<\f{1}{N}.$$ We shall choose $\kp$ so that
$\ds\max_{1\le m\le M}|P_m|<1$. Now we have that
$$ |p_m|\le
\lt|p_m-\sum_{\el=1}^L\kp\,\ap_{m\el}\,x_{\el}\rt|+\kp\lt|\sum_{\el=1}^L\ap_{m\el}\,x_{\el}\rt|
<\f{1}{N}+\kp\sum_{\el=1}^L |\ap_{m\el}|\cdot
|x_{\el}|\le\f{1}{N}+\kp\,L\cdot A\cdot N^{M/L}. $$ Setting
$\ds\f{1}{N}+\kp\,L\cdot A\cdot N^{M/L}=1$, we obtain that
$\ds\kp=\lt(1-\f{1}{N}\rt)\f{1}{L\cdot A\cdot N^{M/L}}$ if $A>0$;
otherwise, it is trivial. Thus we conclude that for $N\ge 2$,
$$\lt|\sum_{\el}^L\ap_{m\el}\,x_{\el}\rt|<\f{1}{\kp N}=\f{L\cdot
A\cdot N^{M/L -1}}{1-\ds\f{1}{N}}\le 2\,L\cdot A\cdot N^{M/L-1}.$$

If $N=1$, then $|x_{\el}|=1$, and hence
$\ds\lt|\sum_{\el=1}^L\ap_{m\el}\,x_{\el}\rt|\le\sum_{\el=1}^L
|\ap_{m\el}|\le L\cdot A.$ \qed

\proclaim{Corollary 1.15 [Carl Ludwig Siegel's Lemma]}

\noindent If the coefficients $a_{m\el}$ of the $M$ linear forms
$\ds\sum_{\el=1}^L a_{m\el}\,x_{\el},\,m=1,2,\cdots,M$, are
integers and $M<L$ , then there is
$(x_1,\cdots,x_L)\in\BZ^L\s\{0\}$ so that $$\sum_{\el=1}^L
a_{1\el}\,x_{\el}=0, \cdots, \sum_{\el=1}^L a_{M\el}\,x_{\el}=0.$$
Moreover, if $A=\ds\max_{1\le\el\le L}\max_{1\le m\le
M}|a_{m\el}|$, then $\ds
|x_{\el}|\le\lt(\lt[(2\,L\,A)^{\f{L}{L-M}}\rt]+1\rt)^{M/L}$.
\endproclaim

\noindent{\it Proof.} Set $2\,L\,A\,N^{M/L -1}<1$ in Corollary
1.14; i.e. $N>(2\,L\,A)^{L/(L-M)}$. Then we take
$N=\ds\lt[(2\,L\,A)^{\f{L}{L-M}}\rt]+1$. \qed

\enddocument